\documentclass[11pt]{article}
\usepackage[utf8]{inputenc}
\usepackage[english]{babel}
\usepackage[active]{srcltx}
\usepackage{hyperref}
\usepackage[all]{xy}
\usepackage{amsfonts}
\usepackage{amsthm}
\usepackage{enumerate}
\usepackage{mathrsfs}
\usepackage{multicol}
\usepackage[all]{xy}
\usepackage{amsfonts}
\usepackage{latexsym,amsmath,amssymb}
\usepackage[usenames,dvipsnames]{color}
\usepackage{pstricks}
\usepackage[absolute,overlay]{textpos}
\usepackage{graphics,wrapfig,times}
\usepackage{xfrac}

\def\tabaddress#1{{\small\it\begin{tabular}[t]{c}#1 \\[1.2ex]\end{tabular}}}
\def\UPCMAT{Department of Mathematics.\\
Universidad Polit\'ecnica de Catalu\~na - Barcelona Tech. \\
   Ed. C-3, Campus Norte.
   C/ Jordi Girona 1. 08034 Barcelona. Spain.}

\newtheorem{teor}{Theorem}
\newtheorem{prop}{Proposition}
\newtheorem{corol}{Corolary}

\newtheorem{definition}{Definition}

\def\beq{\begin{equation}}
\def\eeq{\end{equation}}
\def\bea{\begin{eqnarray}}
\def\eea{\end{eqnarray}}
\def\beann{\begin{eqnarray*}}
\def\eeann{\end{eqnarray*}}
\def\beasn{\begin{sneqnarray}}
\def\eeasn{\end{sneqnarray}}
\def\ben{\begin{enumerate}}
\def\een{\end{enumerate}}
\def\bit{\begin{itemize}}
\def\eit{\end{itemize}}
\def\dst{\(\displaystyle}
\def\proof{( {\sl Proof} )\quad}
\def\derpar#1#2{\frac{\partial{#1}}{\partial{#2}}}

\def\map#1{\mathrel{\mathop{\to}\limits^{#1}}}

\def\feble#1{\mathrel{\mathop\simeq\limits_{#1}}}

\def\qed{\ifvmode\removelastskip\fi
{\unskip\nobreak\hfil\penalty50\hbox{}\nobreak\hfil
\hbox{\vrule height1.2ex width1.2ex}\parfillskip=0pt
\finalhyphendemerits=0 \par\smallskip}}

\def\vf{{\mathfrak X}}
\def\df{{\mit\Omega}}

\def\d{{\rm d}}

\def\Real{{\mathbb R}}

\def\Ker{\mathop{\rm Ker}\nolimits}

\def\inn{\mathop{i}\nolimits}

\def\Tan{{\rm T}}

\def\Lie{\mathop{\rm L}\nolimits}
\def\Cinfty{{\rm C}^\infty}
\catcode`@=12

\title{\sc Connections and jet fields}
\author{\sc Arturo Echeverr\'ia-Enr\'iquez,
   \\
{\sc Miguel C. Mu\~noz-Lecanda\thanks{{\bf e}-{\it mail}: miguel.carlos.munoz@upc.edu}},
   \\
{\sc Narciso Rom\'an-Roy\thanks{{\bf e}-{\it mail}: narciso.roman@upc.edu}},
   \\
   \tabaddress{\UPCMAT}}

\begin{document}
\maketitle
\thispagestyle{empty}

\begin{abstract}
In this review paper
we discuss the different interpretations of the concept of {\sl connection}
in a fiber bundle and in a jet bundle, and its properties,
We relate it with first and second-order systems of partial differential
equations (PDE's) and multivector fields. 
As particular cases we analyze the concepts of linear connections
and connections in a manifold, and their properties and characteristics.
\end{abstract}

 \bigskip
\begin{center}

\vbox{AMS s.\,c.\,(2010): 53B05, 53B15, 53C05, 55R10, 58A20, 58A30.}\null
\end{center}

\clearpage

\tableofcontents

\clearpage

\section{Introduction}

This review paper recovers the contains of several talks given
in an interdisciplinary seminar on themes of Theoretical Physics
and Applied Mathematics.
The main aim is to introduce the concept, characterizations and
properties of {\sl connections} in fiber bundles.
In particular:
\ben
\item
To show the different but equivalent interpretations
of the idea of {\sl connection} in a fiber bundle,
including its relation with {\sl multivector fields}.
\item
To establish the relation between connections in fiber bundles
and systems of partial differential equations and,
in particular, between connections in jet bundles and second-order partial differential equations.
\item
To analyze the characteristics of some particular kinds
of connections: {\sl linear connections} and, as a special case,
{\sl connections on a manifold};
as well as other concepts and properties related to them.
\een

This is a brief review on all these subjects
and all the material presented here is standard an can be found in many books and disertations.
For more information on these and other related topics we address, for instance,
to the references \cite{EMR-98,GHV-72,Hu-66,Ou-73,Sa-89}.

All the manifolds are real, second countable and $\Cinfty$. The maps and the structures are $\Cinfty$.  Sum over repeated indices is understood.

\section{Connections and jet fields in fiber bundles}

In this section, we present the basic elements concerning to first-order jet 
bundles and the theory of jet fields and connections
in fiber bundles (see \cite{Sa-89} for details).

\subsection{First-order jet bundles}

Let $M$ be a differentiable manifold and $\pi\colon E\longrightarrow M$
a differentiable fiber bundle with typical fiber $F$. We denote by 
$\Gamma (M, E)$ or $\Gamma (\pi)$ the set of global sections 
of $\pi$. In the same way, if $U\subset M$ is an open set,
let $\Gamma_U(\pi )$ be the set of local sections of $\pi$
defined on $U$. Let $\dim M=m$ and $\dim E=n$.

For every $y\in E$, the fibers of $J^1\pi$ are denoted $J^1_y\pi$
and their elements by $\bar y$.
If $\phi \colon U\rightarrow E$ is a representative of 
$\bar y \in J^1_y\pi$, we write $\phi \in \bar y$
or $\bar y=\Tan_{\pi (y)}\phi$.
In addition,
the map $\bar\pi^1 = \pi \circ \pi^1\colon J^1\pi \longrightarrow M$
defines another structure of differentiable bundle.
it can be proved \cite{Ou-73} that
$\pi^1 \colon J^1\pi\longrightarrow E$
is an affine bundle modelled on the vector bundle
${\bf E}=\pi^*\Tan^*M \otimes_E{\rm V}(\pi)$
(This notation denotes the tensor product of two vector bundles over $E$).
Therefore, the rank  of $\pi^1\colon J^1\pi\longrightarrow E$ is $mn$. 

We denote by ${\rm V}(\pi)$ and ${\rm V}(\pi^1 )$ the vertical bundles associated with $\pi$ and $\pi^1$ respectively;
that is ${\rm V}(\pi)=\Ker\Tan\pi$ and ${\rm V}(\pi^1)=\Ker\Tan\pi^1$.
We denote by $\vf^{{\rm V}(\pi)}(E)$ and $\vf^{{\rm V}(\pi^1 )}(J^1\pi)$
the corresponding set of sections ; that is, the vertical vector fields.
In the same way we denote by $\vf(E)$ (resp. $\vf(J^1\pi)$) 
the set of vector fields in $E$ (resp. of $J^1\pi$) 
and by $\df^k(E)$ (resp. $\df^k(J^1\pi)$)  the set of differential forms 
of degree $k$ in $E$ (resp. in $J^1\pi$) .

Sections of $\pi$ can be lifted to $J^1\pi$ in the following way:
let $\phi \colon U\subset M \to E$ be a local section of $\pi$,
for every $x\in U$, the section $\phi$ defines  an element of
$J^1\pi$: the equivalence class of $\phi$ in $x$, which is denoted
$(j^1\phi )(x)$.
Therefore we can define a local section $j^1\phi $ of $\bar\pi^1$
and a map $j^1\colon \Gamma_U(\pi)\to \Gamma_U(\bar\pi^1)$
as follows
$$
\begin{array}{ccccc}
j^1\phi &\colon &U & \longrightarrow & J^1\pi \\
& & x & \longmapsto     & (j^1 \phi )(x)
\end{array} 
\quad ; \quad
\begin{array}{ccccc}
j^1&\colon &\Gamma_U(\pi) & \longrightarrow & \Gamma_U(\bar\pi^1)  \\
& & \phi         & \longmapsto     & j^1(\phi )\equiv j^1\phi 
\end{array}  \ .
$$
The section $j^1\phi$ is called the {\rm canonical lifting} or the 
{\sl canonical prolongation} of $\phi$ to $J^1\pi$.
A section of $\bar\pi^1$ which is the canonical extension of a section 
of $\pi$ is called a {\sl holonomic section}.

Let $x^\mu$, $\mu = 1,\ldots,m$, be a local system in $M$
and $y^i$, $i= 1,\ldots,n$, a local system in the fibers;
that is, $\{ x^\mu ,y^i\}$ is a coordinate system adapted to the bundle.
In these coordinates, a local section 
$\phi\colon U \rightarrow E$ is writen as
$\phi (x)=(x^\mu,\phi^i(x))$, that is, $\phi (x)$
is given by functions $y^i=\phi^i(x)$. 
These local systems ${x^\mu}$, ${y^i}$ allows us to construct a  
local system $(x^\mu,y^i,y^i_\mu)$ in $J^1\pi$, where 
$y^i_\mu$ are defined as follows: if $\bar y\in J^1\pi$,
with $\pi^1 (\bar y)=y$ and $\pi(y)=x$, 
let $\phi\colon U\rightarrow E$, $y^i=\phi ^i$, be a representative  
of $\bar y$, then 
$$
y^i_\mu (\bar y)=\left(\derpar{\phi^i}{x^\mu}\right)_x \ .
$$
These coordinate systems are called {\sl natural local systems} in $J^1\pi$.
In them we have
$$
j^1\phi (x)=\left(x^\mu (x),\phi ^i (x),\derpar{\phi^i}{x^\mu}(x)\right) \ .
$$

\subsection{Connections in fiber bundles and jet fields}

In order to set the main definition,
first we prove the following statement:

\begin{teor}
\label{equivtheor}
Let $\pi\colon E\to M$ be a fiber bundle and
$\pi^1\colon J^1\pi\to E$ the corresponding first-order jet bundle.
The following elements can be canonically constructed one from the other:
\ben
\item
A $\pi$-semibasic $1$-form $\nabla$ on $E$
with values in $\Tan E$; that is, an element of
$\Gamma (E,\pi^*\Tan^*M)\otimes\Gamma (E,\Tan E)$,
such that $\nabla^*\alpha =\alpha$, for every
$\pi^1$-semibasic form $\alpha\in\df^1(E)$.
\item
A subbundle ${\rm H}(\nabla)$ of $\Tan E$ such that
\beq
\Tan E={\rm V}(\pi )\oplus{\rm H}(\nabla) \ .
\label{split}
\eeq
\item
A (global) section of $\pi^1\colon J^1\pi\to E$;
that is, a mapping $\Psi\colon E\to J^1\pi$
such that $\pi^1\circ\Psi ={\rm Id}_E$.
\een
\end{teor}
\proof
\quad (1 $\Rightarrow$ 2)\quad
First, observe that 
$\nabla\colon\vf(E)\to\vf(E)$ is a $\Cinfty (E)$-map
which vanishes when it acts on the vertical vector fields.
Its transposed map is $\nabla^*\colon\df^1(E)\to\df^1(E)$,
which is defined as usually by $\nabla^*\beta :=\beta\circ\nabla$,
for every $\beta\in\df^1(E)$. Moreover, as $\nabla$ is
$\pi$-semibasic, so is $\nabla^*\beta$,
then $\nabla^*(\nabla^*\beta )=\nabla^*\beta$
and hence $\nabla\circ\nabla =\nabla$.
Therefore, $\nabla$ and $\nabla^*$ are projection operators
in $\vf(E)$ and $\df^1(E)$ respectively.
So we have the splittings
$$
\vf(E) = {\rm Im}\nabla \oplus{\rm Ker}\nabla
\quad ; \quad
\df^1(E) = {\rm Im}\nabla^* \oplus{\rm Ker}\nabla^* \ .
$$
Now, if ${\cal S}$ is a submodule of $\vf (E)$, the {\sl incident} or
{\sl annihilator} of ${\cal S}$ is defined as the set of $1$-forms
${\cal S}':=\{\alpha\in\df^1(E)\ | \ \alpha (X)=0 \ , \ \forall X\in {\cal S}\}$.
Therefore we have the natural identifications
\beq
({\rm Im}\nabla )'={\rm Ker}\nabla^* \quad ; \quad
({\rm Ker}\nabla )'= {\rm Im}\nabla^* \ .
\label{iden}
\eeq
Taking this into account, for every $y\in E$, the map
$\nabla_y\colon\Tan_yE\to\Tan_yE$ induces the splittings
\beq
\Tan_yE = {\rm Im}\nabla_y \oplus{\rm Ker}\nabla_y
\quad ; \quad
\Tan^*_yE = {\rm Im}\nabla^*_y \oplus{\rm Ker}\nabla^*_y \ .
\label{split1}
\eeq
Next we must prove that ${\rm V}_y(\pi )={\rm Ker}\nabla_y$.
But ${\rm V}_y(\pi )\subseteq{\rm Ker}\nabla_y$, and
${\rm Im}\nabla^*_y$ is the set of $\pi$-semibasic forms at $y\in E$,
then we have
${\rm V}_y(\pi )={\rm Ker}\nabla_y$ and hence
$$
{\rm Ker}\nabla =\Gamma (E,{\rm V}(\pi ))\equiv\vf^{{\rm V}(\pi)}(E)
\quad ; \quad
{\rm Im}\nabla^* =\Gamma (E,\pi^*\Tan^*M) \ .
$$
So we define
$$
{\rm H}(\nabla):=
\bigcup_{y\in E}\{ \nabla_y(u) \ \vert\ u\in\Tan_yE\} \ .
$$
As a consequence of this, the first splitting in (\ref{split1}) leads to
\beq
\Tan E={\rm H}(\nabla)\oplus{\rm V}(\pi ) \ ,
\label{split2}
\eeq
and it allows us to introduce the projections
\beq
{\mathfrak h}\colon\Tan E\longrightarrow{\rm H}(\nabla)\quad ; \quad
{\mathfrak v}\colon\Tan E\longrightarrow{\rm V}(\pi ) \ ,
\label{projections}
\eeq
whose transposed maps
$$
{\mathfrak h}^*\colon{\rm H}^*(\nabla )\longrightarrow\Tan^*E
\quad ; \quad
{\mathfrak v}^*\colon{\rm V}^*(\pi )\longrightarrow\Tan^*E
$$
are injections which lead to the spliting
\beq
\Tan^*E={\rm H}^*(\nabla)\oplus{\rm V}^*(\pi ) \ ;
\label{split3}
\eeq
then, taking into account the second equality of
(\ref{iden}) and (\ref{split1}), in a natural way we have the identifications
${\rm H}(\nabla)'$ with ${\rm V}^*(\pi )$ and
${\rm V}(\pi )'$ with ${\rm H}^*(\nabla)$.

\quad\quad (2 $\Rightarrow$ 1)\quad
Given the subbundle ${\rm H}(\nabla)$ and the splitting
$\Tan E={\rm V}(\pi )\oplus{\rm H}(\nabla)$,
the projections ${\mathfrak h}$ and ${\mathfrak v}$ induce the corresponding
projection operators ${\mathfrak H}$ and ${\mathfrak V}$ in $\vf(E)$
 and the splitting
$X={\mathfrak H}(X)+{\mathfrak V}(X)$, for every $X\in\vf(E)$.
Then we can define the map
$$
\begin{array}{ccccc}
\nabla&\colon&\vf(E)&\longrightarrow&\vf(E)
\\
& & X & \mapsto & {\mathfrak H}(X)
\end{array} \ ,
$$
which is a $\Cinfty (E)$-morphism and satisfies
trivially the following properties:
\ben
\item
$\nabla$ vanishes on the vertical vector fields and
therefore $\nabla\in\Gamma (E,\pi^*\Tan^*M)\otimes\vf(E)$.
\item
$\nabla\circ\nabla=\nabla$,
since $\nabla$ is a projection.
\item
if $\alpha\in\Gamma (E,\pi^*\Tan^*M)$
and $X\in\vf(E)$ we have
$$
(\nabla^*\alpha )X=\alpha (\nabla(X))=
\alpha (h(X))=\alpha ({\mathfrak H}(X)+{\mathfrak V}(X))=\alpha (X)
$$
because $\alpha$ is $\bar\pi^1$-semibasic. 
Therefore $\nabla^*\alpha =\alpha$.
\een
\quad\quad (2 $\Rightarrow$ 3)\quad
Suppose that $\Tan E$ splits as $\Tan E={\rm H}(\nabla)\oplus{\rm V}(\pi)$.
Then, there is a natural way of constructing a section of
$\pi^1\colon J^1\pi\to E$.
In fact, consider $y\in E$ with $\pi(y)=x$, we have
$\Tan_yE={\rm H}_y(E)\oplus{\rm V}_y(\pi)$
and $\Tan_y\pi\vert_{{\rm H}_y(E)}$
is an isomorphism between ${\rm H}_y(E)$ and $\Tan_xM$.
Let $\phi_y\colon U\to E$ be a local section defined in a neigbourhood of
$x$, such that
$$
\phi_y(x)=y \quad  , \quad
\Tan_x\phi_y=(\Tan_y\pi\vert_{{\rm H}_y(E)})^{-1} \ ,
$$
then we have a section
$$
\begin{array}{ccccc}
\Psi&\colon&E&\longrightarrow&J^1\pi
\\
& &y&\mapsto&(j^1\phi_y)(\pi(y))
\end{array} \ ,
$$
which is differentiable because the splitting
$\Tan_yE={\rm H}_y(E)\oplus{\rm V}_y(\pi)$
depends differentiabily on $y$.

(3 $\Rightarrow$ 2)\quad
Let $\Psi\colon E\to J^1\pi$ be a section
and $\bar y\in J^1\pi$, with $\bar y\map{\pi^1} y\map{\pi} x$.
Observe that $\Psi (y)\in J^1\pi$ is an equivalence class
of sections $\phi\colon M\to E$, with $\phi (x)=y$,
but the subspace ${\rm Im}\Tan_x\phi$
does not depend on the representative $\phi$, provided it is in this class.
Then, for every $\bar y\in J^1\pi$ and being
$\phi$ a representative of $\bar y=\Psi (y)$, we define
$$
{\rm H}_y(\nabla):={\rm Im}\Tan_x\phi
\quad {\rm and}\quad
{\rm H}(\nabla):=\bigcup_{y\in E}{\rm H}_y(\nabla) \ .
$$
\qed

\begin{definition}
A {\rm connection} in the bundle $\pi\colon E\to M$
is any of the equivalent elements of Theorem \ref{equivtheor}.
Then, a (global) section $\Psi\colon E\to J^1\pi$ is said to be a
{\rm jet field} in the bundle $\pi^1\colon J^1\pi\to E$.
The $\pi$-semibasic form $\nabla$ is called the {\rm connection form}
or {\rm Ehresmann connection}.
The subbundle ${\rm H}(\nabla)$ is called the {\rm horizontal subbundle} of
$\Tan E$ associated with the connection and
the sections of ${\rm H}(\nabla)$  are the {\rm horizontal vector fields}.
It is also denoted ${\cal D}(\Psi)$ and is called the
{\rm distribution associated} with $\Psi$.
\end{definition}

A jet field $\Psi\colon E\to J^1E$
(resp. an Ehresmann connection $\nabla$) is said to be {\sl orientable}
if ${\cal D}(\Psi)$ is an orientable distribution on $E$.
If $M$ is orientable, then every connection in $E$ is also orientable.

\noindent{\bf Remarks}:
\bit
\item
If $\bar y\in J^1_yE$, with $x=\pi(y)$, and $\phi\colon U\to E$  
is a representative of $\bar{y}$, we have the split
$$
\Tan_yE={\rm Im}\ \Tan_x\phi\oplus{\rm V}_y(\pi) \ .
$$
Hence the sections of $\pi^1$ are identified with connections 
in the bundle $\pi\colon E\longrightarrow M$,
since they induce a horizontal subbundle of $\Tan E$.
Observe that it is reasonable to write ${\rm Im}\ \bar y$ for an 
element $\bar y\in J^1\pi$.
\item
Any global section of an affine bundle  can be identified with its associated vector bundle.
In particular:
\begin{itemize}
\item
Let $\pi\colon E\longrightarrow M$ be a trivial bundle; that is
$E=M\times F$. A section of $\pi^1$ can be chosen in the following way:
denoting by $\pi_1\colon M\times F\to M$ and $\pi_2\colon M\times F\to F$
the canonical projections,
for a given $y_o\in M\times F$, $y_o=(x_o,v_o)=(\pi_1(y_o),\pi_2(y_o))$,
we define the section 
$\phi_{y_o}(x)=(x,\pi_2(y_o))$, for every $x\in M$.
From a section of $\pi$ we construct another one of $\pi^1$ as follows:
$$
z(y):=(j^1\phi_y)(\pi_1(y))\quad ; \quad y\in E \ ,
$$
which is taken as the zero section of $\pi^1$.
In this case, $J^1\pi$ is a vector bundle over $E$.
\item
If $\pi\colon E\longrightarrow M$ is a vector bundle with
typical fiber $F$, let $\phi :M\to E$ be
the zero section of $\pi$ and $j^1\phi\colon M\to J^1\pi$
its canonical lifting. We construct the zero section of
$\bar\pi^1$ in the following way: 
$$
z(y):=(j^1\phi)(\pi(y))\quad ; \quad y\in E \ ;
$$
thereby, in this case $\bar\pi\colon J^1\pi\longrightarrow M$
is a vector bundle.
\end{itemize}
\eit

\subsection{Local expressions and properties}
\label{locexp}

Let $(x^\mu ,y^i)$ be a local system of coordinates in an open set
$U\subset E$. The most general local expression of a
semibasic $1$-form on $E$ with values in $\Tan E$ is
$$
\nabla=
f_\mu\d x^\mu\otimes\left( g^\nu\derpar{}{x^\nu}+h^i\derpar{}{y^i}\right) \ .
$$
As $\nabla^*$ is the identity on semibasic forms, it follows that
$\nabla^*\d x^\mu =\d x^\mu$,
so the local expression of the connection form $\nabla$ is
$$
\nabla=
\d x^\mu\otimes\left(\derpar{}{x^\mu}+{\mit\Gamma}_\mu^i\derpar{}{y^i}\right)
\ ,
$$
where ${\mit\Gamma}^i_\mu\in\Cinfty (U)$.
In this system the jet field $\Psi$ is expressed as
$$
\Psi =(x^\mu ,y^i,{\mit\Gamma}_\rho^i(x^\mu ,y^i)) \ .
$$

Let $\phi$ be a representative of $\Psi (y)$ with $\phi =(x^\mu ,f^i(x^\mu ))$.
Therefore $\phi (x)=y$, $\Tan_x\phi=\Psi (y)$ and we have
$$
y=\phi (x)=(x^\mu ,f^i(x^\mu ))=(x^\mu ,y^i) \ .
$$
The matrix of $\Tan_x\phi$ is
$\displaystyle
\left(\begin{array}{c} {\rm  Id} \\ \left(\displaystyle\derpar{f^i}{x^\mu}\right)_x\end{array}\right)$,
therefore
\dst\derpar{f^i}{x^\nu}\Big\vert_x={\mit\Gamma}_\nu^i(x^\mu ,y^i)\) .
Now, taking \dst\left\{\derpar{}{x^\mu}\Big\vert_x\right\}\) as a basis of
$\Tan_xM$, we obtain
$$
{\rm Im}\Tan_x\phi =
\left\{ (\Tan_x\phi )\left(\derpar{}{x^\mu}\Big\vert_x\right)\right\} =
\left\{\derpar{}{x^\mu}\Big\vert_y+
{\mit\Gamma}_\mu^i(y)\derpar{}{y^i}\Big\vert_y\right\} \ ,
$$
hence, ${\rm H}(\nabla)$ is locally spaned by
\beq
\left\{\derpar{}{x^\mu}+{\mit\Gamma}_\mu^i(y)\derpar{}{y^i}\right\} \ .
\label{lochor}
\eeq

As final remarks, notice that the
splitting (\ref{split2}) induces a further one
$$
\vf(E) = {\rm Im}\nabla \oplus\Gamma (E,{\rm V}(\pi )) \ ;
$$
so every vector field $X\in\vf (E)$ splits
into its {\sl horizontal} and {\sl vertical} components:
$$
X=X^{\rm H}+X^{\rm V}=\nabla (X)+(X-\nabla (X)) \ ;
$$
that is, ${\mathfrak H}\equiv\nabla$ and ${\mathfrak V}\equiv{\rm Id}-\nabla$.
Locally, this splitting is given by
$$
X = f^\mu\derpar{}{x^\mu}+g^i\derpar{}{y^i} =
f^\mu\left(\derpar{}{x^\mu}+{\mit\Gamma}_\mu^i\derpar{}{y^i}\right)+
\left( g^i-f^\mu{\mit\Gamma}_\mu^i\right)\derpar{}{y^i} \ ,
$$
since \dst\derpar{}{x^\mu}+{\mit\Gamma}_\mu^i\derpar{}{y^i}\) and
\dst\derpar{}{y^i}\) generate locally $\Gamma (E,{\rm H}(\nabla))$ and
$\Gamma (E,{\rm V}(\pi ))$, respectivelly.
Observe that, if $X$ is an horizontal vector field,
then $\nabla (X)=X$.

In an analogous way the splitting (\ref{split3}) induces the following one
$$
\df^1(E) = \Gamma (E,\pi^*\Tan^*M)\oplus{\rm Ker}\nabla^* =
\Gamma (E,\pi^*\Tan^*M)\oplus ({\rm Im}\nabla )' \ ;
$$
then, for every $\alpha\in\df^1 (E)$, we have
$$
\alpha=\alpha^{\rm H}+\alpha^{\rm B}=
\nabla^*\alpha+(\alpha -\nabla^*\alpha ) \ ,
$$
whose local expression is
\beq
\alpha = F_\mu\d x^\mu+G_i\d y^i =
(F_\mu +G_i{\mit\Gamma}_\mu^i)\d x^\mu+
G_i(\d y^i-{\mit\Gamma}_\mu^i\d x^\mu ) \ ,
\label{vfsplit}
\eeq
since $\d x^\mu$ and $\d y^i-{\mit\Gamma}_\mu^i\d x^\mu$
generate locally $\Gamma (E,{\rm H}^*(\nabla))$ and
$\Gamma (E,{\rm V}^*(\pi ))$, respectivelly.

As a final remark, we analyze the structure
of the set of connections in $\pi\colon E\to M$.
Then, let $\nabla_1,\nabla_2$ be two connection forms.
The condition $\nabla_1^*\alpha =\nabla_2^*\alpha=0$,
for every semibasic $1$-form $\alpha$, means that
$(\nabla_1-\nabla_2)^*\alpha=0$; that is
$\nabla_1-\nabla_2\in\Gamma (E,\pi^*\Tan^*M)\otimes_E\Gamma (E,{\rm V}(\pi ))$.

However, let $\nabla$ be a connection on $\pi\colon E\to M$ and
$\gamma\in\Gamma (E,\pi^*\Tan^*M)\otimes_E\Gamma (E,{\rm V}(\pi ))$,
then $\nabla +\gamma$ is another connection form. So we have:

\begin{prop}
The set of connection forms on $\pi\colon E\to M$ is an
affine ``space'' over the module of semibasic differential $1$-forms on $E$
with values in ${\rm V}(\pi )$.
\end{prop}

In a local canonical system, if
\dst\nabla =\d x^\mu\otimes
\left(\derpar{}{x^\mu}+{\mit\Gamma}^i_\mu\derpar{}{y^i}\right)\)
and \dst\gamma=\gamma^i_\mu\d x^\mu\otimes\derpar{}{y^i}\) \ ,
then
$$
\nabla +\gamma =\d x^\mu\otimes
\left(\derpar{}{x^\mu}+({\mit\Gamma}^i_\mu +\gamma^i_\mu )\derpar{}{y^i}\right) \ .
$$

\subsection{Integrability of jet fields and connections.
First-order partial differential equations}

\begin{definition}
The {\rm curvature} of a connection $\nabla$ is
a $(2,1)$-tensor field in $E$
which is defined as follows: for every $Z_1,Z_2\in\vf (E)$,
$$
{\cal R}(Z_1,Z_2):=
({\rm Id}-\nabla )([\nabla (Z_1),\nabla (Z_2)])=
\inn ([\nabla (Z_1),\nabla (Z_2)])({\rm Id}-\nabla ) \ .
$$
\label{curva}
\end{definition}

Using the coordinate expressions of the
connection form $\nabla$ or the jet field $\Psi$,
a simple calculation leads to
$$
{\cal R} = \frac{1}{2}
\left(\derpar{{\mit\Gamma}_\nu^j}{x^\mu}-
\derpar{{\mit\Gamma}_\mu^j}{x^\nu}+
{\mit\Gamma}_\mu^i\derpar{{\mit\Gamma}_\nu^j}{y^i}-
{\mit\Gamma}_\nu^i\derpar{{\mit\Gamma}_\mu^j}{y^i}\right)
(\d x^\mu\wedge\d x^\nu )\otimes\derpar{}{y^j} \ .
$$

\begin{definition}
Let $\Psi\colon E\to J^1\pi$ be a jet field associated with a connection $\nabla$.
\ben
\item
A section $\phi\colon M\to E$ is an
{\rm integral section} of $\Psi$ (resp. of $\nabla$) if
$\Psi\circ\phi =j^1\phi $.
\item
$\Psi$ is an {\rm integrable jet field}
(resp. $\nabla$ is an {\rm integrable connection})  if it admits
integral sections.
\een
\end{definition}

One may readily check that, if $(x^\mu ,y^i,y_\mu^i)$ is a natural local
system in $J^1\pi$ and, in this system,
$\Psi =(x^\mu ,y^i,{\mit\Gamma}_\rho^i(x^\mu ,y^i))$ and
$\phi =(x^\mu ,f^i(x^\nu ))$,
then $\phi$ is an integral section of $\Psi$ if, and only if,
$\phi$ is a solution of the following system of
partial differential equations
\beq
\derpar{f^i}{x^\mu}={\mit\Gamma}_\mu^i\circ\phi \ .
\label{condin}
\eeq

The integrable jet fields and connections can be characterized as follows:

\begin{prop}
The following assertions on a jet field $\Psi$ are equivalent:
\ben
\item
The jet field $\Psi$ is integrable.
\item
The curvature of the connection form $\nabla$
associated with $\Psi$ is zero.
\item
${\cal D}(\Psi )$ is an involutive distribution.
\een
\end{prop}
( {\sl Proof} )\quad
(1 $\Leftrightarrow$ 2)\quad
Notice that if $\phi$ is an integral section of $\Psi$, then the distribution
${\cal D}(\Psi )$ is tangent to the image of $\phi$, and conversely.

\quad\quad (2 $\Leftrightarrow$ 3)\quad
From the definition (\ref{curva}) we obtain that, if ${\cal R}=0$, then
$$
\nabla([\nabla (Z_1),\nabla (Z_2)])=[\nabla (Z_1),\nabla (Z_2)] \ ,
$$
hence, the horizontal distribution ${\cal D}(\Psi )$ is involutive.
Conversely, if ${\cal D}(\Psi )$ is involutive,
as $\nabla$ is the identity on ${\cal D}(\Psi )$,
the last equation follows, and then ${\cal R}=0$.
\qed

\noindent{\bf Remark}:
According to this proposition,
from the local expression of ${\cal R}$ we obtain the local
integrability conditions of the equations (\ref{condin}).

\subsection{Connections and multivector fields}

\begin{definition}
A {\rm $k$-multivector field} in $E$ is a section of 
$\Lambda^m(\Tan E)=\overbrace{\Tan E\wedge\ldots\wedge\Tan E}^k$
or, what is equivalent, a skew-symmetric contravariant 
tensor of order $k$ in $E$. The set of $k$-multivector fields 
in $E$ is denoted $\vf^k(E)$.

A $k$-multivector field $\mathbf{X}\in\vf^k(E)$ is said to be {\rm locally decomposable} if,
for every $y\in E$, there is an open neighbourhood  $U_y\subset E$
and $X_1,\ldots ,X_k\in\vf (U_y)$ such that $\mathbf{X}\vert_{U_y}=X_1\wedge\ldots\wedge X_k$.
\end{definition}

\begin{definition}
If $\Omega\in\df^r(E)$ and ${\bf X}\in\vf^k(E)$,
the {\rm contraction} between ${\bf X}$ and $\Omega$ is defined as
the natural contraction between tensor fields. In particular,
for locally decomposable multivector fields,
$$
 \inn({\bf X})\Omega\mid_{U}=\inn(X_1\wedge\ldots\wedge X_k)\Omega 
=\inn (X_1)\ldots\inn (X_k)\Omega \ .
$$
\end{definition}

\begin{definition}
 A $k$-multivector field $\mathbf{X}\in\mathfrak{X}^k(E)$ is 
 \textsl{$\pi$-transverse} if, for every $\beta\in\Omega^k(M)$ such that
$\beta (\pi(y))\not= 0$, at every point
$y\in E$, we have that  $(\inn(\mathbf{X})(\pi^*\beta))_y\not= 0$. 
\end{definition}

Let ${\cal D}$ be a $k$-dimensional distribution in $E$;
that is, a $k$-dimensional subbundle of $\Tan E$.
Obviously sections of $\Lambda^m{\cal D}$ are $k$-multivector fields in $E$.
The existence of a non-vanishing global section of $\Lambda^k{\cal D}$
is equivalent to the orientability of ${\cal D}$. Therefore:

\begin{definition}
A non-vanishing $k$-multivector field ${\bf X}\in\vf^k(E)$ and
a $k$-dimensional distribution ${\cal D}\subset\Tan E$
are {\rm locally associated} if there exists a connected open set
$U\subseteq E$ such that ${\bf X}\vert_U$ is a section of $\Lambda^m{\cal D}\vert_U$.
\end{definition}

As a consequence of this definition,
if ${\bf X},{\bf X}'\in\vf^k(E)$ are non-vanishing multivector fields
locally associated with the same distribution ${\cal D}$, 
on the same connected open set $U$, then there exists a
non-vanishing function $f\in\Cinfty (U)$ such that
${\bf X}'\feble{U}f{\bf X}$. This fact defines an equivalence relation in the
set of non-vanishing $k$-multivector fields in $E$, whose equivalence classes
are denoted by $\{ {\bf X}\}_U$. Then:

\begin{teor}
There is a bijective correspondence between the set of $k$-dimensional
orientable distributions ${\cal D}$ in $\Tan E$ and the set of the
equivalence classes $\{ {\bf X}\}_E$ of non-vanishing, locally decomposable
$k$-multivector fields in $E$.
\label{bijcor}
\end{teor}
\proof
Let $\omega\in\df^k(E)$ be an orientation form for ${\cal D}$.
If $y\in E$ there exists an open neighbourhood $U_y\subset E$
and $X_1,\ldots ,X_k\in\vf (U_p)$,
with $\inn(X_1\wedge\ldots\wedge X_k)\omega >0$, such that
${\cal D}\vert_{U_y}=\langle X_1,\ldots ,X_k\rangle$.
Then $X_1\wedge\ldots\wedge X_k$ is a representative of a class
of $k$-multivector fields associated with ${\cal D}$ in $U_y$.
But the family $\{ U_y\ ;\ y\in E\}$ is a covering of $E$;
let $\{ U_\alpha\ ;\ \alpha\in A\}$ be
a locally finite refinement and $\{ \rho_\alpha\ ;\ \alpha\in A\}$
a subordinate partition of unity.
If $X^\alpha_1,\ldots ,X^\alpha_k$ is a local basis of ${\cal D}$ in $U_\alpha$,
with $\inn(X^\alpha_1\wedge\ldots\wedge X^\alpha_k)\omega >0$,
then \dst {\bf X}=\sum_\alpha\rho_\alpha X^\alpha_1\wedge\ldots\wedge X^\alpha_k\)
is a global representative of the class of non-vanishing
$k$-multivector fields associated with ${\cal D}$ in $E$.

The converse is immediate since, if
${\bf X}\vert_U=X^1_1\wedge\ldots\wedge X^1_k=X^2_1\wedge\ldots\wedge X^2_k$,
for different sets $\{ X^1_1,\ldots ,X^1_k\}$ and $\{X^2_1,\ldots ,X^2_k\}$,
then
$\langle X^1_1,\ldots ,X^1_k\rangle =\langle X^2_1,\ldots ,X^2_k\rangle$.
\qed

\begin{definition}
A $k$-multivector field
$\mathbf{X}\in\vf^k(E)$ is {\rm integrable} if its associated distribution 
${\cal D}({\bf X})$ is integrable. 
Then the integral submanifolds of $\mathbf{X}\in\vf^k(E)$ are the integral submanifolds of ${\cal D}({\bf X})$.
\end{definition}

If ${\bf X}\in\vf^k(E)$ is locally decomposable,
then ${\bf X}$ is $\pi$-transverse if, and only if,
$\Tan_y\pi({\cal D}({\bf X}))=\Tan_{\pi (y)}M$, for every $y\in E$.
(Remember that ${\cal D}({\bf X})$ is the $k$-distribution associated to ${\bf X}$).

\begin{teor}
\ben
\item
Let ${\bf X}\in\vf^m(E)$ be integrable.
Then ${\bf X}$ is $\pi$-transverse if, and only if,
its integral manifolds are local sections of $\pi\colon E\to M$.
\item
${\bf X}\in\vf^m(E)$ is integrable and $\pi$-transverse if, and only if, for every point $y\in E$,
there exists a local section $\phi\colon U\subset M\to E$
such that $\phi(\pi(y))=y$, and a non-vanishing
function $f\in\Cinfty (E)$ such that
$\Lambda^m\Tan\phi =f{\bf X}\circ\phi\circ\Lambda^k\tau_U$.
\een
\label{insecmvf1}
\end{teor}
\proof
\ben
\item
Consider $y\in E$, with $\pi(y)=x$. In a neighbourhood of $y$ there exist
$X_1,\ldots ,X_k\in\vf (E)$ such that $X_1,\ldots ,X_m$ span ${\cal D}({\bf X})$ and
${\bf X}=X_1\wedge\ldots\wedge X_m$. But, as $X$ is $\pi$-transverse,
$(\inn (X)(\pi^*\omega))_y\not= 0$,
for every $\omega\in\df^m(M)$ with $\omega (x)\not= 0$.
Thus, taking into account the second comment above,
${\cal D}({\bf X})$ is a $\pi$-transverse distribution and
$X_\mu\not\in\vf^{{\rm V}(\pi)}(E)$ at any point, for every $X_\mu$.
Now, let $S\hookrightarrow E$ be the integral manifold of ${\cal D}({\bf X})$
passing through $y$, then
$\Tan_yS=\langle(X_1)_y,\ldots ,(X_m)_y\rangle$.
As a consequence of all of this, and again taking into account
the second comment above, for every point $y\in S$,
$\Tan_y\pi({\cal D}({\bf X}))=\Tan_{\pi (y)}M$,
then $\pi\vert_S$ is a local diffeomorphism
and $S$ is a local section of $\pi$.
The converse is obvious.
\item
If ${\bf X}$ is integrable and $\pi$-transverse, then by theorem \ref{insecmvf1},
for every $y\in E$, with $\pi(y)=x$,
there is an integral local section $\phi\colon U\subset M\to E$
of ${\bf X}$ at $y$ such that $\phi (x)=y$.
Then, as a consequence of the definition of integrability,
 ${\bf X}_y$ spans $\Lambda^m\Tan_y\phi$,
and hence the relation in the statement holds.

Conversely, if the relation holds, then
${\rm Im}\,\phi$ is an integral manifold of ${\bf X}$ at $y$,
then ${\bf X}$ is integrable and, as $\phi$ is a section of $\pi$,
${\bf X}$ is necessarily $\pi$-transverse.
\qed
\een

In this case, if $\phi\colon U\subset M\to E$
is a local section with $\phi (x)=y$ and $\phi (U)$ is
the integral manifold of ${\bf X}$ through $y$,
then $\Tan_y({\rm Im}\,\phi)$ is ${\cal D}_y({\bf X})$.

\begin{teor}
There is a bijective correspondence between the set of
orientable jet fields $\Psi\colon E\to J^1\pi$
(that is, the set of orientable Ehresmann connection $\nabla$ in $\pi\colon E\to M$)
and the set of the equivalence classes
of locally decomposable and $\pi$-transverse $m$multivector fields
$\{ {\bf X}\}\subset\vf^m(E)$.
They are characterized by the fact that
${\cal D}(\Psi)={\cal D}({\bf X})$.

In addition, the orientable jet field $\Psi$ is integrable if, and only if, so is ${\bf X}$,
for every ${\bf X}\in\{ {\bf X}\}$.
\label{vfmvf1}
\end{teor}
\proof
If $\Psi$ is an orientable jet field in $J^1\pi$,
let ${\cal D}(\Psi)$ its horizontal distribution.
Then, taking ${\cal D}\equiv{\cal D}(\Psi)$, we construct
$\{ {\bf X}\}$ by applying theorem \ref{bijcor}
and, since the distribution ${\cal D}(\Psi)$ is $\pi$-transverse,
the result follows immediately.
The proof of the converse statement is similar.

Moreover, $\Psi$ is integrable if, and only if,
${\cal D}(\Psi)={\cal D}({\bf X})$ is also. 
Therefore it follows that ${\bf X}$ is also integrable, for ${\bf X}\in\{ {\bf X}\}$,
and conversely.
\qed

Reminding the local expressions of Section  (\ref{locexp}),
we have that the local expression for
a particular representative multivector field ${\bf X}$ of the class
$\{ {\bf X}\}\subset\vf^k(E)$ associated with a jet field $\Psi$ (or a connection form $\nabla$) is
$$
{\bf X}\equiv\bigwedge_{\mu=1}^k X_\mu=\bigwedge_{\mu=1}^k
\left(\derpar{}{x^\mu}+{\mit\Gamma}_\mu^A\derpar{}{y^i}\right) \ .
$$
Then, $\phi =(x^\mu ,f^i (x^\nu ))$ is an integral section of ${\bf X}$ 
if, and only if,
$\phi$ is a solution of the system of partial differential equations
(\ref{condin}).

\section{Connections and jet fields in jet bundles}

The geometrical framework for treating with
systems of second order partial differential equations are the
 jet bundles $J^1\pi$ and $J^1J^1\pi$.
Next we analize this topic.

\subsection{Connections in $J^1\pi$ and jet fields in $J^1J^1\pi$}

Consider the bundle $\bar\pi^1\colon J^1\pi\to M$.
The jet bundle $J^1J^1\pi$ is obtained
by defining an equivalence relation on the local sections
of $\bar\pi^1$. Hence, the elements of $J^1J^1\pi$ are
equivalence classes of these local sections and
$J^1J^1\pi$ is an affine bundle over $J^1\pi$,
modelled on the vector bundle
$\bar\pi^{1^*}\Tan^*M \otimes_{J^1\pi}{\rm V}(\pi^1 )$.   
So, we have the commutative diagram
\beq
\begin{array}{ccccc}
J^1J^1\pi&\put(-10,3){\vector(1,0){25}}\put(-2,8){\mbox{$\pi^1_1$}}
&J^1\pi&\put(-10,3){\vector(1,0){25}}\put(-2,8){\mbox{$\pi^1$}}&E \\
&
\begin{picture}(20,20)(0,0)
\put(-2,0){\mbox{$\bar\pi^1_1$}}
\put(0,20){\vector(1,-1){25}}
\end{picture}
&\put(0,20){\vector(0,-1){25}}\put(3,7){\mbox{$\bar\pi^1$}}&
\begin{picture}(20,20)(0,0)
\put(13,0){\mbox{$\pi$}}
\put(25,20){\vector(-1,-1){25}}
\end{picture}
&
\\
& & M & &
\end{array}
\label{diag}
\eeq

Let ${\cal Y}\colon J^1\pi\to J^1J^1\pi$ be a jet field in $J^1J^1\pi$.
As we know ${\cal Y}$ induces a connection form
$\nabla$ and a horizontal $n+1$-subbundle ${\rm H}(\nabla)$ such that
$$
\Tan J^1\pi ={\rm V}(\pi^1 )\oplus{\rm Im}\,\nabla =
{\rm V}(\pi^1 )\oplus{\rm H}(\nabla) \ ,
$$
where ${\rm H}_{\bar y}={\rm Im}\Tan_{\bar\pi^1 (\bar y)}\psi$,
for $\bar y\in J^1\pi$ 
and $\psi\colon M\to J^1\pi$ a representative of ${\cal Y}(\bar y)$.
We denote by ${\cal D}({\cal Y})$ the $\Cinfty (J^1\pi)$-module of
sections of ${\rm H}(\nabla)$.

If $(x^\mu ,y^i,y_\mu^i,z_\nu^i,z_{\nu\mu}^i)$ is a natural
system of coordinates in $J^1J^1\pi$, we have the following local expressions
for these elements
\beann
{\cal Y}&=&(x^\mu ,y^i,y_\mu^i,F_\nu^i(x^\rho ,y^j,y_\rho^j),
G_{\nu\mu}^i(x^\rho ,y^j,y_\rho^j)) \ ,
\\
{\rm H}(J^1\pi)&=&
\left\{\derpar{}{x^\mu}+F_\mu^i(\bar y)\derpar{}{y^i}+
G_{\nu\mu}^i(\bar y)\derpar{}{y_\nu^i}\right\} \ ,
\\
\nabla&=&\d x^\mu\otimes
\left(\derpar{}{x^\mu}+F_\mu^i\derpar{}{y^i}+
G_{\nu\mu}^i\derpar{}{y_\nu^i}\right) \ .
\eeann

From these local expressions
we obtain that a representative $m$-multivector field ${\bf X}$ of the class
$\{ {\bf X}\}\subset\vf^m(J^1\pi)$ associated with the jet field ${\cal Y}$, 
has the local expression
$$
{\bf X}\equiv\bigwedge_{\mu=1}^m X_\mu=
\bigwedge_{\mu=1}^m \left(\derpar{}{x^\mu}+F_\mu^i\derpar{}{y^i}+
G_{\mu\nu}^i\derpar{}{y_\nu^i}\right) \ .
$$

If ${\cal Y}\colon J^1\pi\to J^1J^1\pi$ is a jet field then
a section $\psi \colon M\to J^1\pi$ is an
{\sl integral section} of ${\cal Y}$ if
${\cal Y}\circ \psi =j^1\psi $.
${\cal Y}$ is an {\sl integrable jet field} if it admits
integral sections.
One may readily check that,
in a natural system of coordinates in $J^1J^1\pi$,
if $\psi =(x^\mu ,f^\mu (x^\nu ),g^i_\mu (x^\nu )$, then
it is an integral section of ${\cal Y}$ if, and only if,
$\psi $ is a solution of the following system of differential equations
$$
\derpar{f^i}{x^\mu}=F_\mu^i\circ\psi
\qquad \derpar{g_\rho^i}{x^\mu}=G_{\rho\mu}^i\circ\psi \ .
$$
Remember that if $\psi$ is an integral section of ${\cal Y}$,
then the distribution
${\cal D}({\cal Y})$ is tangent to the image of $\psi$ and conversely.
Hence, ${\cal Y}$ is integrable if, and only if,
${\cal D}({\cal Y})$ is an involutive distribution
or, what is equivalent, if, and only if,
the curvature of $\nabla$ is zero; that is, in coordinates
\beann
0 &=&
\left(\derpar{F_\eta^j}{x^\mu}+F_\mu^i\derpar{F_\eta^j}{y^i}+
G_{\gamma\mu}^i\derpar{F_\eta^j}{y_\gamma^i}-
\derpar{F_\mu^j}{x^\eta}-F_\eta^i\derpar{F_\mu^j}{y^i}-
G_{\rho\eta}^i\derpar{F_\mu^j}{y_\rho^i}\right)
(\d x^\mu\wedge\d x^\eta )\otimes\derpar{}{y^j}+
\\ & &
\left(\derpar{G_{\rho\eta}^j}{x^\mu}+F_\mu^i\derpar{G_{\rho\eta}^j}{y^i}+
G_{\gamma\mu}^i\derpar{G_{\rho\eta}^j}{y_\gamma^i}-
\derpar{G_{\rho\mu}^j}{x^\eta}-F_\eta^i\derpar{G_{\rho\mu}^j}{y^i}-
G_{\gamma\eta}^i\derpar{G_{\rho\mu}^j}{v_\gamma^i}\right)
(\d x^\mu\wedge\d x^\eta )\otimes\derpar{}{y_\rho^j} \ .
\eeann

\subsection{The SOPDE condition. Holonomic jet fields and connections}

The idea of this Section is to characterize the integrable jet fields
in $J^1J^1\pi$ such that their integral sections
are canonical prolongations of sections of the projection $\pi$.

It is well known that there are two natural projections from
$\Tan\Tan Q$ to $\Tan Q$. In the same way, if we
consider the diagram (\ref{diag}), we see that
there is another natural projection from $J^1J^1\pi$ to $J^1\pi$.
Let ${\bf y}\in J^1J^1\pi$ with
\dst{\bf y}\stackrel{\pi^1_1}{\mapsto}\bar y\stackrel{\pi^1}{\mapsto}
y\stackrel{\pi}{\mapsto}x\) ,
and $\psi\colon M\to J^1\pi$ a representative of ${\bf y}$, that is,
${\bf y} =\Tan_x\psi$. Consider now the section
$\phi =\pi^1\circ\psi\colon M\to E$, then
$j^1\phi (x)\in J^1\pi$ and we have:

\begin{prop}
The following projection is a differentiable map:
$$
\begin{array}{ccccc}
j^1\pi^1&\colon&J^1J^1\pi&\longrightarrow&J^1\pi
\\
& &{\bf y}&\mapsto&j^1(\pi^1\circ\psi )(\bar\pi^1_1({\bf y})) 
\end{array} \ .
$$
\end{prop}
( {\sl Proof} )\quad
Let $(x^\mu ,y^i,y_\mu^i,z_\nu^i,z_{\nu\mu}^i)$ be a natural coordinate
system in a neigbourhood of ${\bf y}_0\in J^1J^1\pi$ and
${\bf y}_0=({x_0}^\mu,{y_0}^i,{y_0}_\mu^i,{z_0}_\nu^i,{z_0}_{\nu\mu}^i)$.
We have $\pi^1_1({\bf y}_0)=({x_0}^\mu,{y_0}^i,{y_0}_\mu^i)$.
Let $\psi\colon M\to J^1\pi$ be a representative of ${\bf y}_0$;
locally $\psi=(x^\mu,f^i(x^\nu ),g_\mu^i(x^\nu ))$ with
$$
f^i({x_0})={y_0}^i \quad ,\quad g_\mu^i({x_0})={y_0}^i_\mu
\quad ;\quad
\derpar{f^i}{x^\nu}(x_0)={z_0}_\nu^i \quad ,\quad
\derpar{g_\mu^i}{x^\nu}(x_0)={z_0}^i_{\mu\nu}
$$
then $\phi =\pi^1\circ\psi =(x^\mu ,f^i(x^\nu ))$, and
$$
j^1\phi =\left( x^\mu ,f^i(x^\rho),\derpar{f^i}{x^\nu}(x^\rho)\right)
$$
Hence $j^1\pi^1({\bf y}_0)=({x_0}^\mu,{y_0}^i,{z_0}^i_\nu)$
and the result follows.
\qed

\noindent{\bf Remark}:
Observe that $j^1\pi^1$ and $\pi^1_1$ exchange the coordinates
$y^i_\mu$ and $z^i_\mu$.

\begin{corol}
If $\psi\colon M\to J^1\pi$ is a section of $\bar\pi^1$,
then $j^1\pi^1\circ j^1\psi =j^1(\pi^1\circ\psi )$.
\end{corol}
\proof
In a coordinate system $(x^\mu ,y^i,y_\mu^i,z_\nu^i,z_{\nu\mu}^i)$,
we have $\psi =(x^\mu ,f^i(x^\rho),g_\nu^i(x^\rho ))$ and
$$
j^1\psi =\left( x^\mu ,f^i(x^\rho),g_\nu^i(x^\rho),
\derpar{f^i}{x^\nu}(x^\rho),
\derpar{g_\mu^i}{x^\nu}(x^\rho)\right) \ ,
$$
but
\dst j^1\pi^1\circ j^1\psi =\left( x^\mu ,f^i(x^\rho),
\derpar{f^i}{x^\rho}(x^\nu)\right) =j^1(\pi^1\circ\psi )\) .
\qed

\begin{definition}
A jet field ${\cal Y}\colon J^1\pi\to J^1J^1\pi$ is  a
{\rm Second Order Partial Differential Equation (SOPDE)}
(or also that verifies the {\rm SOPDE condition}) if it
is a section of the projection $j^1\pi^1$ or, what is equivalent,
$$
j^1\pi^1\circ {\cal Y}={\rm Id}_{J^1\pi} \ .
$$
\end{definition}

Now, we are going to characterize SOPDE integrable jet fields. First we define:

\begin{definition}
A section $\psi\colon M\to J^1E$ of $\bar\pi^1$ is said to be {\rm  holonomic} if it is the canonical lifting of a section $\phi\colon M\to E$ of $\pi$; that is,
$\psi=j^1\phi$.
\end{definition}

\begin{prop}
Let ${\cal Y}\colon J^1\pi\to J^1J^1\pi$ be an integrable jet field.
The necessary and sufficient condition for ${\cal Y}$ to be a SOPDE is that
its integral sections are holonomic.
\end{prop}
\proof
\quad ($\Longrightarrow$)\quad
If ${\cal Y}$ is a SOPDE then  $j^1\pi^1\circ {\cal Y}={\rm Id}_{J^1\pi}$.
Let $\psi\colon M\to J^1\pi$ be an integral section of ${\cal Y}$; that is,
${\cal Y}\circ\psi =j^1\psi$, then
$$
\psi =j^1\pi^1\circ{\cal Y}\circ\psi =
j^1\pi^1\circ j^1\psi =j^1(\pi^1\circ\psi )\equiv j^1\phi \ ,
$$
and $\psi$ is a canonical prolongation.

\quad\quad ($\Longleftarrow$)\quad
Now, let ${\cal Y}$ be an integrable jet field whose integral sections
are canonical prolongations. Take $\bar y\in J^1\pi$ and
$\phi\colon M\to E$ a section such that
$j^1\phi\colon M\to J^1\pi$ is an integral section of ${\cal Y}$
with $j^1\phi (\bar\pi^1(\bar y))=\bar y$. We have
\beann
(j^1\pi^1\circ{\cal Y})(\bar y)&=&
(j^1\pi^1\circ{\cal Y})(j^1\phi (\bar\pi^1(\bar y)))=
(j^1\pi^1\circ{\cal Y}\circ j^1\phi )(\bar\pi^1(\bar y))=
(j^1\pi^1\circ{\cal Y}\circ\psi )(\bar\pi^1(\bar y))
\\ &=&
(j^1\pi^1\circ j^1\psi)(\bar\pi^1(\bar y))=
j^1(\pi^1\circ\psi)(\bar\pi^1(\bar y))=
j^1\phi (\bar\pi^1(\bar y))=\bar y \ ,
\eeann
and ${\cal Y}$ is a SOPDE.
\qed

\noindent{\bf Remarks}:
\bit
\item
In coordinates, the condition $j^1\pi^1\circ {\cal Y}={\rm Id}_{J^1\pi}$
is expressed as follows: the jet field
${\cal Y}=(x^\mu ,y^i,y_\mu^i,F_\nu^i,G_{\nu\mu}^i)$
is a SOPDE if, and only if, $F_\nu^i=y_\nu^i$.
\item
If ${\cal Y}=(x^\mu ,y^i,y_\mu^i,y_\nu^i,G_{\nu\mu}^i)$
is a SOPDE, then 
\dst j^1\phi =\left( x^\mu ,f^i,\derpar{f^i}{x^\mu}\right)\)
is an integral section of ${\cal Y}$ if, and only if,
$\phi$ is the solution of the system of (second order) PDE's 
\beq
G_{\nu\mu}^i\left(x^\rho ,f^j,\derpar{f^j}{x^\gamma}\right) =
\frac{\partial^2f^i}{\partial x^\nu\partial x^\mu} \ ,
\label{sisteq2}
\eeq
which justifies the terminology.
Sometimes, SOPDE jet fields which are not integrable are also called {\sl semi-holonomic jet fields}.
\eit

If ${\cal Y}$ is a SOPDE, then a representative 
$m$-multivector field ${\bf X}$ of the class
$\{ {\bf X}\}\subset\vf^m(J^1\pi)$ associated with ${\cal Y}$, 
has the local expression
$$
{\bf X}\equiv
\bigwedge_{\mu=1}^m \left(\derpar{}{x^\mu}+y_\mu^i\derpar{}{y^i}+
G_{\mu\nu}^i\derpar{}{y_\nu^i}\right) \ ,
$$
and, if ${\bf X}$ is integrable,
its (holonomic) integral sections are the solutions of (\ref{sisteq2}).
 It is important to remark that,
since the integrability of a class of multivector fields is equivalent to demanding that
the curvature ${\cal R}$ of the connection associated with this class
vanishes everywhere; the system (\ref{sisteq2}) has solution
if, and only if, the following additional system of equations holds
\beann
0 &=& G^j_{\nu\mu}-G^j_{\mu\nu} \ ,
\\
0 &=&  \derpar{G_{\nu\rho}^j}{x^\mu}+v_\mu^i\derpar{G_{\nu\rho}^j}{y^i}+
G_{\mu\gamma}^i\derpar{G_{\nu\rho}^j}{y_\gamma^i}-
\derpar{G_{\mu\rho}^j}{x^\nu}-v_\nu^i\derpar{G_{\mu\rho}^j}{y^i}-
G_{\nu\gamma}^i\derpar{G_{\mu\rho}^j}{y_\gamma^i} \ .
\eeann

\section{Connections in a vector bundle}

As special cases of connections in fiber bundles,
we study connections in a vector bundle,
in particular {\sl linear connections} and
the related notion of {\sl connections on a manifold}.

\subsection{Structures in a vector bundle}

Linear connections are a particular type of connections, which
can be defined only on vector bundles.
In order to give their different characterizations,
first we need to introduce some previous concepts.

First, remember that, if ${\cal E}$ is a vector space
and $u\in{\cal E}$, there is a natural identification between
${\cal E}$ and $\Tan_u{\cal E}$, which is given by
$$
\begin{array}{ccc}
{\cal E} & \longrightarrow & \Tan_u{\cal E} \\
w & \mapsto & {\rm D}_w(u)
\end{array} \ ,
$$
where ${\rm D}_w(u)$ is the directional derivative with respect
the vector $w$ at the point $u$; that is, for every
differentiable function $f:{\cal E}\to \Real$,
$$
{\rm D}_w(u)f:=\lim_{t\to 0}\frac{1}{t}(f(u+tw)-f(u)) \ .
$$
So, if $z^1,\ldots ,z^n$ are coordinates in ${\cal E}$
and $w\equiv (\lambda^1,\ldots ,\lambda^n)$, then
\dst{\rm D}_w(u)=\lambda^i\derpar{}{z^i}\Big\vert_u\) ,
and the identification is immediate.

In this way, if $\pi\colon E\to M$ is a vector bundle,
for every section $\phi\colon M\to E$ and $p\in M$,
we denote by $\sharp_{\phi (p)}\colon{\rm V}_{\phi (p)}(\pi )\to E_p$
the natural identification between the vector space
$E_p$ and ${\rm V}_{\phi (p)}(\pi )=\Tan_{\phi (p)}E_p$.

\begin{definition}
Let $\phi (p)\equiv (p,u)\in E$, with
$u\in{\rm V}_{\phi (p)}(\pi )=\Tan_uE_p$.
The map
\beann
E & \longrightarrow & \Tan E={\rm V}(\pi ) \\
(p,u)&\mapsto&((p,u),{\rm D}_u(u))
\eeann
defines a $\pi$-vertical vector field $\Delta\in\vf (E)$
which is called the {\rm Liouville vector field} of the
vector bundle $\pi\colon E\to M$.
\end{definition}

In a natural set of coordinates $(x^\mu ,y^i)$ of $E$,
the local expression of $\Delta$ is \dst\Delta =y^i\derpar{}{y^i}\) .
Then, its integral curves are the solution of the system of
differential equations
$$
\frac{\d x^\mu}{\d t}=0 \quad , \quad \frac{\d y^i}{\d t}=y^i \ ;
$$
that is, $x^\mu(t)=A^\mu$ and $y^i(t)=B^ie^t$ ($A^\mu ,B^i\in\Real$).
Thus, $\Delta$ generates dilatations along the fibres of the
vector bundle.

\noindent{\bf Remark}:
Observe that the flow of $\Delta$ is the map
$$
\begin{array}{ccccc}
\Phi&\colon&\Real\times E&\longrightarrow&E
\\
& &(t,(x,y))&\mapsto&(x,e^ty)
\end{array} \ ;
$$
which, by a suitable reparametrization, can be equivalently defined as
$$
\begin{array}{ccccc}
\Phi&\colon&\Real^+\times E&\longrightarrow&E
\\
& &(t,(x,y))&\mapsto&(x,ty)
\end{array} \ .
$$
This map allows to define the following one: for every $t\in \Real^+$,
the map $\Phi_t\colon E\to E$ is given by $\Phi_t(y)=ty$, for every $y\in E$.

If a connection is given in a vector bundle,
then we can establish the following:

\begin{definition}
Let $\nabla$ be a connection on a vector bundle $\pi\colon E\to M$.
The {\rm covariant derivative} induced by $\nabla$ is the map
$\tilde\nabla\colon\Gamma (\pi)\to\Gamma (\pi)\otimes\df^1(M)$
defined as follows: if $\phi\colon M\to E$ is a section,
$Z\in\vf (M)$ and $p\in M$,
$$
(\tilde\nabla\phi )(p;Z):= \sharp_{\phi (p)}({\mathfrak v}(\Tan_p\phi (Z_p)))=
\sharp_{\phi (p)}[\Tan_p\phi (Z_p)-\nabla (\phi(p);\Tan_p\phi (Z_p))] \ .
$$
(It is usual to write $(\tilde\nabla_Z\phi )(p)$ instead of
$(\tilde\nabla\phi )(p;Z)$).
\end{definition}

\subsection{Linear connections}

Now we can prove the following equivalence:

\begin{teor}
Let $\pi\colon E\to M$ be a vector bundle.
Consider a connection given by the equivalent elements
$\nabla$, ${\rm H}(\nabla)$ or $\Psi$.
Then the following conditions are equivalent:
\ben
\item
The connection form is invariant under the Liouville vector field:
$$
\Lie (\Delta )\nabla =0 \ ,
$$
or, what is the same thing,
the {\sl vertical projection operator} ${\mathfrak V}\equiv {\rm Id}-\nabla$ is invariant by $\Delta$:
$$
\Lie (\Delta ){\mathfrak V}=0 \ .
$$
\item
The Liouville vector field preserves the horizontal subbundle;
that is, for every $t\in\Real^+$ and every $y\in E$, we have
$$
\Tan_y\Phi_t({\rm H}_y(\nabla))={\rm H}_{\Phi_t(y)}(\nabla) \ .
$$
\item
If $(x^\mu ,y^i)$ is a bundle system of coordinates
in the vector bundle $\pi\colon E\to M$, then
the functions ${\mit\Gamma}^i_\mu$ which characterize the connection
are linear on the fibers and their expressions are
${\mit\Gamma}^i_\mu =\pi^*(-{\mit\Gamma}^i_{j\mu})y^j$,
where ${\mit\Gamma}^i_{j\mu}\in\Cinfty (M)$ are the {\rm Christoffel symbols}
of the linear connection.
\item
The jet field $\Psi$ is a vector bundle morfism.
(Notice that if $E\to M$ is a vector bundle, so is $J^1\pi\to M$).
\item
For every $f\in\Cinfty (M)$ and every section $\phi\colon M\to E$,
$$
\tilde\nabla (f\phi )=\d f\otimes\phi +f\tilde\nabla\phi \ .
$$
\een
\label{linconn}
\end{teor}
\proof
\quad (1 $\Rightarrow$ 2)\quad
If $\Lie (\Delta )\nabla =0$ then $\Phi_{t_*}\nabla =\nabla$,
for every $t\in\Real^+$.
If $u\in{\rm H}_y(E)$, then there exists $v\in\Tan_yE$ such that
$\nabla (v)=u$ and we have
$$
\Tan_y\Phi_t(u)=\Tan_y\Phi_t(\nabla (v))=(\Phi_{t_*}\nabla )(\Phi_{t_*}v)=
\nabla (\Phi_{t_*}v) \ ,
$$
which implies that $\Tan_y\Phi_t(u)\in{\rm H}_{\Phi_t(y)}(\nabla)$
and, since $\Phi_t$ is a diffeomorphism the result follows.

\quad\quad (2 $\Rightarrow$ 1)\quad
If $\Tan_y\Phi_t({\rm H}_y(\nabla))={\rm H}_{\Phi_t(y)}(\nabla)$,
for every $t\in\Real^+$, then the splitting
$\Tan_yE={\rm H}_y(E)+{\rm V}_y(\pi )$ implies that
$$
\Tan_{\Phi_t(y)}E={\rm H}_{\Phi_t(y)}(\nabla)+{\rm V}_{\Phi_t(y)}(\pi )=
\Tan_y\Phi_t({\rm H}_y(\nabla))+{\rm V}_{\Phi_t(y)}(\pi )=
\Tan_y\Phi_t{\rm H}_y(\nabla)+\Tan_y\Phi_t({\rm V}_y(\pi )) \ ;
$$
that is, if $u\in\Tan_yE$, writing $u={\mathfrak h}(u)+{\mathfrak v}(u):=\nabla (u)+{\mathfrak v}(u)$, we have
$$
\Tan_y\Phi_t(u)=\nabla(\Tan_y\Phi_t(u))+{\mathfrak v}(\Tan_y\Phi_t(u))=
\Tan_y\Phi_t(\nabla (u))+\Tan_y\Phi_t({\mathfrak v}(u)) \ ,
$$
and hence $\nabla\circ\Tan_y\Phi_t=\Tan_y\Phi_t\circ\nabla$,
so $\nabla$ is invariant by $\Phi_t$ and, therefore,
$\Lie (\Delta )\nabla =0$.

\quad\quad (1 $\Leftrightarrow$ 3)\quad
Locally we have
\dst\nabla =
\d x^\mu\otimes\left(\derpar{}{x^\mu}+\Gamma^i_\mu\derpar{}{y^i}\right)\) ,
then $\Lie (\Delta )\nabla =0$ implies that the functions
$\Gamma^i_\mu$ are homogeneous of degree $1$ ({\sl Euler's theorem})
and, as they are differentiable at the origin, then they are
also linear in the variables $y^i$. So
${\mit\Gamma}^i_\mu =\pi^*(-{\mit\Gamma}^i_{j\mu})y^j$,
where ${\mit\Gamma}^i_{j\mu}\in\Cinfty(M)$ 
(they are functions of the coordinates $x^\mu$).

The converse is immediate.

\quad\quad (3 $\Leftrightarrow$ 4)\quad
Taking into account that the local expression of $\Psi$ is
$\Psi (x^\mu ,y^i)=(x^\mu ,y^i,\Gamma^i_\mu(x,y))$,
the assertion is immediate since $\Psi$ is a vector bundle morphism
if, and only if, the functions $\Gamma^i_\mu$ are linear of
the coordinates $y^i$.

\quad\quad (3 $\Leftrightarrow$ 5)\quad
On the one hand we have that
\beann
(\d f\otimes\phi +f\nabla\phi )(p;Z)&=&
Z_p(f)\phi(p)+f(p)(\tilde\nabla\phi)(p;Z)
\\ &=&
Z_p(f)\phi(p)+f(p)(\sharp_{\phi (p)}\circ {\mathfrak v}\circ\Tan_p\phi )Z_p
\\ &=&
Z_p(f)\phi(p)+f(p)\sharp_{\phi (p)}(\Tan_p\phi )Z_p -
(\nabla (\phi (p));\Tan_p\phi Z_p)) \ ,
\eeann
and, on the other hand,
$$
\tilde\nabla (f\phi )(p;Z)=
\sharp_{f(p)\phi (p)}[\Tan_p(f\phi )Z-\nabla ((f\phi )(p);\Tan_p(f\phi )Z_p] \ .
$$
In a local natural system of coordinates in the vector bundle,
if \dst Z=g^\mu\derpar{}{x^\mu}\) and $\phi =(x^\mu ,\phi^i )$,
then
\beann
\Tan_p(f\phi )Z_p&=&
g^\mu (p)\derpar{}{x^\mu}\Big\vert_{(f\phi )(p)}+
f(p)\derpar{\phi^i}{x^\mu}\Big\vert_p
g^\mu (p)\derpar{}{y^i}\Big\vert_{(f\phi )(p)}
\\ & &
+g^\mu (p)\derpar{f}{x^\mu}\Big\vert_p
\phi^i(p)\derpar{}{y^i}\Big\vert_{(f\phi )(p)} \ ,
\\
\tilde\nabla (f\phi)(p;Z)&=&
\left( p;g^\mu (p)\left( f(p)\derpar{\phi^i}{x^\mu}\Big\vert_p+
\derpar{f}{x^\mu}\Big\vert_p\phi^i(p)-
{\mit\Gamma}^i_\mu (f(p)\phi (p))\right)\right) \ ,
\\
(\d f\otimes\phi +f\tilde\nabla\phi )(p;Z)&=&
\left( p;g^\mu (p)\left( f(p)\derpar{\phi^i}{x^\mu}\Big\vert_p+
\derpar{f}{x^\mu}\Big\vert_p\phi^i(p)-
f(p){\mit\Gamma}^i_\mu (\phi (p))\right)\right) \ ;
\eeann
therefore, the local condition in order that item 5 holds is
$$
{\mit\Gamma}^i_\mu (f(p)\phi (p))=f(p){\mit\Gamma}^i_\mu (\phi (p)) \ ,
$$
for every $f\in\Cinfty (M)$. This proves the assertion.
\qed

Then we define

\begin{definition}
A connection in the vector bundle $\pi\colon E\to M$
is a {\rm linear connection} if
the equivalent conditions in Theorem \ref{linconn} hold.
\end{definition}

Finally, we can state the following concept:

\begin{definition}
Let $\nabla$ be a linear connection on a vector bundle $\pi\colon E\to M$,
and $\phi\colon M\to E$ a section. Then $\nabla\phi$ is called the
{\rm covariant differential} of $\phi$.
It is a map $\nabla\phi\colon\vf (M)\to\Gamma (\pi)$,
which is an element of $\df^1(M)\otimes_M\Gamma (\pi)$.
\end{definition}

In natural coordinates the local expression of $\nabla\phi$ is
$$
\nabla\phi =
\d x^\mu\otimes\left(x^\rho ,\derpar{\phi^i}{x^\mu}+
{\mit\Gamma}^i_{j\mu}\phi^j\right) \ .
$$


In relation to the structure of the set of linear connections,
if $\nabla_1, \nabla_2$ are linear connections, then
$\nabla_1-\nabla_2\colon\Gamma (\pi)\to\vf (M)\otimes\Gamma (\pi)$
is a $\Cinfty (M)$-linear map, so that
$$
\nabla_1-\nabla_2\in\Gamma (\pi)^*\otimes\df^1(M)\otimes\Gamma (\pi)=
\df^1(M)\otimes{\rm End}\, \Gamma (\pi)=
\df^1(M)\otimes_M\Gamma (M,{\rm L}_E) \ ,
$$
where ${\rm L}_E$ is the bundle of endomorphisms of $E$.
On the other hand, if
${\Upsilon}\in\df^1(M,{\rm L}_E)=\df^1(M)\otimes\Gamma (M,{\rm L}_E)$
and $\nabla$ is a linear connection, then
$\nabla +{\Upsilon}$ is another linear connection because
$$
(\nabla +{\Upsilon} )(f\phi )=\d f\otimes (\nabla +{\Upsilon} )+
f(\nabla +{\Upsilon} )\phi \ ,
$$
since the action of ${\Upsilon}$ on sections is the following:
writing ${\Upsilon} =\alpha_i\otimes\upsilon^i$, with
$\alpha_i\in\df^1(M)$ and $\upsilon^i\in\Gamma(M,{\rm L}_E)$,
and taking $Z\in\vf (M)$, $\phi\colon M\to E$ and $p\in M$, then
$$
{\Upsilon} (p;Z,\phi )=\alpha_i (p;Z)\upsilon^i(\phi (p))
$$
verifying that
$$
{\Upsilon} (p;Z,f\phi )=f(p){\Upsilon} (p;Z,\phi )\ ,
$$
because ${\Upsilon}$ is $\Cinfty (M)$-linear.
Therefore we can state:

\begin{prop}
The set of linear connections on $\pi\colon E\to M$
is an affine ``space'' modelled on the $\Cinfty (M)$-module of $1$-forms on $M$
with values on the bundle of endomorphisms of $E$, ${\rm L}_E$.
\end{prop}

As you can observe, there exists a canonical injection of the module
$\df^1(M)\otimes_M\Gamma (M,{\rm L}_E)$ into
$\Gamma (E,\pi^*\Tan^*M)\otimes_E\Gamma (E,{\rm V}(\pi ))$
defined as follows:
$\alpha\otimes{\Upsilon}\mapsto
\pi^*\alpha\otimes (\sharp^{-1}\circ{\Upsilon})$,
and that this injection is a morphism of $\Cinfty (M)$-modules.

\section{Connections in a manifold}

\subsection{Basic definitions and properties and covariant derivative}

The concept of {\sl connection on a manifold} is closely
related to that of linear connection.
Let $M$ be an $m$-dimensional differentiable manifold.

\begin{definition}
A {\rm connection} on $M$ is
a linear connection in the tangent bundle $\Tan M$.
\end{definition}

In this case, the sections of the bundle are vector fields,
hence we can define:

\begin{definition}
Let $\nabla$ be a connection on the manifold $M$.
\ben
\item
If $X,Y\in\vf (M)$, the map
$$
{\cal T}(X,Y):=\tilde\nabla_XY-\tilde\nabla_YX-[X,Y]
$$
is called the {\rm torsion} associated to this connection.
${\cal T}$ is an antisymmetric tensor field on $M$ with values
on $\Tan M$; so that ${\cal T}\in\df^2(M)\otimes\vf (M)$.
\item
$\nabla$ is a {\rm torsion-free} or {\rm symmetric connection} if
${\cal T}$ is zero.
\een
\end{definition}

Let $(x^\mu ,v^\mu )$ be a natural system of coordinates of $\Tan M$.
If the local expression of $\nabla$ is
$$
\nabla =\d x^\mu\otimes\left(\derpar{}{x^\mu}+
{\mit\Gamma}^\rho_{\mu\nu}v^\nu\derpar{}{x^\rho}\right) \ ,
$$
then, for every vector field in $M$
(a section of the vector bundle $\Tan M\to M$) we have
$$
\tilde\nabla_{\derpar{}{x^\rho}}\left( g^\eta\derpar{}{x^\eta}\right) =
\left(\derpar{g^\eta}{x^\rho}+
{\mit\Gamma}^\eta_{\rho\nu}g^\nu\right)\derpar{}{x^\nu} \ ,
$$
and the local expression of ${\cal T}$ is
$$
{\cal T}=(-{\mit\Gamma}^\mu_{\rho\eta}+{\mit\Gamma}^\mu_{\eta\rho})\,
\d x^\rho\otimes\d x^\eta\otimes\derpar{}{x^\mu} \ ,
$$
According to this we can state the following Proposition
(which justifies the name given to these kinds of connections):

\begin{prop}
The necessary and sufficient condition for the connection $\nabla$ on $M$
to be torsion-free is that the following relation holds
for every system of coordinates,
$$
{\mit\Gamma}^\mu_{\rho\eta}={\mit\Gamma}^\mu_{\eta\rho} \ .
$$
\end{prop}

If $\nabla_1,\nabla_2$ are symmetric connections on $M$,
this implies that, for every $X,Y\in\vf (M)$,
$$
\tilde\nabla_{1_X}Y-\tilde\nabla_{1_Y}X=
\tilde\nabla_{2_X}Y-\tilde\nabla_{2_Y}X \ ,
$$
therefore
$(\tilde\nabla_1-\tilde\nabla_2)_XY=(\tilde\nabla_1-\tilde\nabla_2)_YX$
and hence $\nabla_1-\nabla_2$ is symmetric too; that is,
it is an element of ${\rm S}_2(M)\otimes\vf (M)$,
where ${\rm S}_2(M)$ denotes the module of
symmetric $2$-degree tensor fields on $M$.
(Remember that an element of ${\rm S}_2(M)\otimes\vf (M)$
is a section of the bundle
$\vee^2\Tan^*M\otimes\Tan M\to M$,
where $\vee^2\Tan^*M$ is the symmetric product of $\Tan^*M$
with itself.
If $(x^\mu )$ is a local system of coordinates in $M$, then a local
section of this bundle is locally expressed as
\dst{\mit\Gamma}^\gamma_{\mu\nu}
\d x^\mu\otimes\d x^\nu\otimes\derpar{}{x^\gamma}\) ,
with ${\mit\Gamma}^\gamma_{\mu\nu}={\mit\Gamma}^\gamma_{\nu\mu}$).

Conversely, if $\nabla$ is a symmetric connection on $M$ and
${\mit\Sigma}\in{\rm S}_2(M)\otimes\vf (M)$,
then $\nabla+{\mit\Sigma}$ is another symmetric connection on $M$.
Therefore:

\begin{prop}
The set of symmetric connections on a manifold $M$
is an affine space modelled on the $\Cinfty (M)$-module
${\rm S}_2(M)\otimes_M\vf (M)$.
\end{prop}

Observe that this module is a submodule of
$\df^1(M)\otimes_M\Gamma (M,{\rm L}_E)=
\df^1(M)\otimes_M\df^1(M)\otimes_M\vf (M)$.


\subsection{Covariant derivative along a path and parallel transport}
Here we recall some elementary constructions on a manifold $M$ with a connection $\nabla$. Let $\tau\colon\Tan M\to M$ be the natural projection.

Let $\sigma \colon I=(-\varepsilon,\varepsilon)\subset\Real\to M$ be a smooth curve. A \textbf{vector field along} $\sigma$ is a smooth maping $V:I\to \mathfrak{X}(M)$ such that $ \tau\circ V=\sigma$. The set of all vector fields along $\sigma$ is denoted by $\mathfrak{X}(M;\sigma)$. With the natural operations, $\mathfrak{X}(M;\sigma)$ is module over the algebra of smooth real functions defined on the interval $I\subset\Real$. The elements in $\mathfrak{X}(M;\sigma)$ can be understood as curves in $\Tan M$. If $V\in\mathfrak{X}(M;\sigma)$ we can say that $V$ is a lifting of $\sigma$ from $M$ to $\Tan M$.

As it is known, the covariant derivative $\nabla_XY\in\mathfrak{X}(M)$, for $X,Y\in\mathfrak{X}(M)$, can be extended to a covariant derivative of vector fields along a curve with respect to the tangent vector to the curve, that is $\nabla_{\dot{\sigma}}V$, for $V\in\mathfrak{X}(M;\sigma)$, obtaining another element of $\mathfrak{X}(M;\sigma)$ with the natural properties with respect to the real functions defined on the interval $I\subset\Real$.

An element  $V\in\mathfrak{X}(M;\sigma)$ is said to be \textbf{parallel along the curve} if  $\nabla_{\dot{\sigma}}V=0$. Given $u_p\in \Tan_p$, with $p=\sigma(0)$, there exists only one element $V\in\mathfrak{X}(M;\sigma)$  satisfying the conditions
\begin{enumerate}
  \item $\nabla_{\dot{\sigma}}V=0$.
  \item $V(0)=u_p$. 
\end{enumerate}
With this in mind it is easy to prove that the set of vector fields along a curve $\sigma$, thai is $\mathfrak{X}(M;\sigma)$ is a module over the real functions defined on the same interval than the curve with finite rang equal to the dimension of the manifold $M$. 

\subsection{Horizontal liftings and covariant derivatives}

Let $\nabla$ be a connection in $M$.

\begin{definition}
Let $p\in M$, $u\in\Tan_pM$, and $\mbox{\boldmath $\sigma$} \colon (-\varepsilon,\varepsilon)\subset\Real\to M$
a smooth curve with $\mbox{\boldmath $\sigma$}(0)=p$.
The {\rm horizontal lifting} of $\mbox{\boldmath $\sigma$}$ to 
the point $u_p=(p,u)\in\Tan M$ is the curve
${\cal X}\colon(-\varepsilon,\varepsilon)\subset\Real\to\Tan M$, 
which is solution to the initial value problem given by
$$
(i)\quad \nabla_{\dot{\mbox{\boldmath $\sigma$}}}{\cal X}=0
\quad , \quad
(ii)\quad {\cal X}(0)=u_p \quad ;
$$
that is, ${\cal X}$ is parallel along $\mbox{\boldmath $\sigma$}$ and
coincides with $u_p$ at $p$.
(Observe that ${\cal X}$ is actually a vector field along  $\mbox{\boldmath $\sigma$}$;
i.e., ${\cal X}\in\vf(M;\mbox{\boldmath $\sigma$})$), 
\end{definition}

In a local chart of coordinates $(x^\mu)$ in $M$,
if $\mbox{\boldmath $\sigma$}=(\sigma^1,\ldots,\sigma^m)$ is
the local expression of $\mbox{\boldmath $\sigma$}$ in this chart, and
$$
{\cal X}(t)={\cal X}^\mu(t)\derpar{}{x^\mu}\Big\vert_{\mbox{\boldmath $\sigma$}(t)} \quad , \quad
\dot{\mbox{\boldmath $\sigma$}}(t)=
\dot\sigma^\mu(t)\derpar{}{x^\mu}\Big\vert_{\mbox{\boldmath $\sigma$}(t)} \ ,
$$
taking into account that
$\displaystyle\nabla_{\derpar{}{x^\mu}}\left(\derpar{}{x^\nu}\right)=\Gamma_{\mu\nu}^\rho\derpar{}{x^\rho}$,
then the local expression of $\nabla_{\dot{\mbox{\boldmath $\sigma$}}}{\cal X}=0$ is
$$
\nabla_{\dot{\mbox{\boldmath $\sigma$}}}{\cal X}=
\dot {\cal X}^\rho(t)+\Gamma_{\nu\mu}^\rho(\mbox{\boldmath $\sigma$}(t))\,{\cal X}^\mu(t)\dot\sigma^\nu(t)=0 \ ,
$$
and the unique solution$ {\cal X}^\rho(t)$ gives the horizontal lifting of $\mbox{\boldmath $\sigma$}$.
Thus, in a natural chart of coordinates $(x^\mu,v^\mu)$ in $\Tan M$, the curve
${\cal X}$ has the expression
${\cal X}(t)=(\sigma^\mu(t),{\cal X}^\mu(t))$,
and its tangent vector at every point is
$$
\dot{\cal X}(t)=(\dot\sigma^\mu(t),\dot{\cal X}^\mu(t))=
(\dot\sigma^\mu(t),-\Gamma_{\nu\mu}^\rho(\mbox{\boldmath $\sigma$}(t))\,{\cal X}^\mu(t)\dot\sigma^\nu(t))\ .
$$
In particular,
$$
\dot{\cal X}(0)=
(\dot\sigma^\mu(0),-\Gamma_{\mu\nu}^\rho(p)\,u^\mu\dot\sigma^\nu(0))\equiv
\dot\sigma^\mu(0)\,\derpar{}{x^\rho}\Big\vert_{u_p}-
\Gamma_{\mu\nu}^\rho(p)\,u^\mu\dot\sigma^\nu(0)\,\derpar{}{v^\rho}\Big\vert_{u_p}\ ,
$$
which depends only on the tangent vector to the curve at $p$.

As a consequence, given $v\in\Tan_pM$, if the curve $\mbox{\boldmath $\sigma$}$
is a representative of $v$ (that is, $\mbox{\boldmath $\sigma$}(0)=p, \dot{\mbox{\boldmath$\sigma$}}(0)=v)$,
then $\dot{\cal X}(0)\in\Tan_{u_p}(\Tan M)$ depends only on $v$ and not 
on the selected representative curve.

\begin{definition}
The {\rm horizontal lifting} of the vector $v\in\Tan_pM$ to $u\in\Tan_pM$ is the vector
tangent at $u_p\in\Tan_pM$ to the horizontal lifting of any curve representative of $v$.

The map that implements this operation is denoted
${\rm h}_p^{u_p}\colon\Tan_pM\to\Tan_{u_p}(\Tan M)$.
\end{definition}

Locally, if $\displaystyle v=v^\rho\derpar{}{x^\rho}\Big\vert_{u_p}$ and
$\displaystyle u=u^\rho\derpar{}{x^\rho}\Big\vert_{u_p}$, we have that
\beq
{\rm h}_p^{u_p}(v)=v^\rho\derpar{}{x^\rho}\Big\vert_{u_p}-
\Gamma_{\nu\mu}^\rho(p)\,u^\mu v^\nu\,\derpar{}{v^\rho}\Big\vert_{u_p}\ .
\label{horlift}
\eeq

\begin{prop}
The map ${\rm h}_p^{u_p}\colon\Tan_pM\to\Tan_{u_p}(\Tan M)$ has the following properties:
\ben
\item
It is a linear map.
\item
It is an injective map.
\item
${\rm Im}\,{\rm h}_p^{u_p}$ is an $m$-dimensional vector subspace of $\Tan_{u_p}(TM)$.
\item
$\Tan_{u_p}(TM)={\rm V}_{u_p}(\tau)\oplus {\rm Im}\,{\rm h}_p^{u_p}$.
\een
\end{prop}
\proof
These properties are an immediate consequence of the local expression \eqref{horlift}. In particular:
\ben
\item
It is a consequence of the linearity of all the operations.
\item
It holds because $\Tan_{u_p}\tau\circ{\rm h}_p^{u_p}={\rm Id}_{\Tan_pM}$.
\item
It is a consequence of  the above items (1) and (2).
\item
It is a consequence of  the above items (1), (2) and (3).
\qed
\een

\begin{definition}
${\rm Im}\,{\rm h}_p^{u_p}$ is the {\rm horizontal subspace} in $u\in\Tan_pM$
associated with the connection $\nabla$,
and it is denoted ${\rm H}_{u_p}\tau$.
\end{definition}

A basis for ${\rm H}_{u_p}\tau$ is given by taking a basis
$\displaystyle\left\{\derpar{}{x^\nu}\Big\vert_p \right\}$ in $\Tan_pM$ and
obtaining the corresponding horizontal liftings (from \eqref{horlift}):
$$
{\rm h}_p^{u_p}\left(\derpar{}{x^\nu}\Big\vert_p \right)=
\derpar{}{x^\nu}\Big\vert_{u_p}-
\Gamma_{\nu\mu}^\rho(p)\,u^\mu\,\derpar{}{v^\rho}\Big\vert_{u_p}=
\derpar{}{x^\nu}\Big\vert_{u_p}-
\Gamma_{\nu\mu}^\rho(p)\,v^\mu(u_p)\,\derpar{}{v^\rho}\Big\vert_{u_p}\ .
$$
Observe also that $\Tan_{u_p}\tau\colon{\rm H}_{u_p}\to\Tan_pM$ is an isomorphism.

The expression of ${\rm h}_p^{u_p}$ depends differentialy on
$p,u_p$ and $v_p$; then the subspace ${\rm H}_{u_p}\tau$
depends differentialy on $u_p$ and hence it defines a subbundle of $\Tan(\Tan M)$ of rank $m$, which is denoted
$\displaystyle{\rm H}(\nabla):=\bigcup_{{u_p}\in\Tan M}{\rm H}_{u_p}\tau$,
and is called the {\sl horizontal subbundle} associated with the connection $\nabla$.
Obviously we have that $\Tan(\Tan M)={\rm V}(\tau)\oplus{\rm H}(\nabla)$.

A local basis of ${\rm H}(\nabla)$ is given by the vector fields
$$
\displaystyle \left\{ \derpar{}{x^\nu}-
\Gamma_{\nu\mu}^\rho\,v^\mu\,\derpar{}{v^\rho}\right\}\ ;
$$
hence, comparing this local expression with \eqref{lochor},
and taking $\Gamma_\nu^\rho=-\Gamma_{\mu\nu}^\rho\,v^\mu$
(as stated in the item 3 of Theorem \ref{linconn}),
we see that this horizontal subspace is just the one introduced in \eqref{split2}.

In this way, we have the two projections \eqref{projections}:
$$
{\mathfrak h}\colon\Tan(\Tan M)\longrightarrow{\rm H}(\nabla)\quad ; \quad
{\mathfrak v}\colon\Tan(\Tan M)\longrightarrow{\rm V}(\tau ) \ ,
$$
whose extension to vector fields ${\mathfrak H}$ and ${\mathfrak V}$ 
have the local expressions are obtained from \eqref{vfsplit} and, in this case, 
if $\displaystyle\left\{\derpar{}{x^\nu},\derpar{}{v^\nu}\right\}$ is a local basis for $\vf(\Tan M)$, are explicitly given by
\beann
{\mathfrak H}\left(\derpar{}{x^\nu}\right)=
\derpar{}{x^\nu}-\Gamma_{\nu\mu}^\rho\,v^\mu\,\derpar{}{v^\rho}
\quad &,& \quad
{\mathfrak V}\left(\derpar{}{x^\nu}\right)=
\Gamma_{\nu\mu}^\rho\,v^\mu\,\derpar{}{v^\rho} \quad ;
\\
{\mathfrak H}\left(\derpar{}{v^\nu}\right)=0
\quad &,& \quad
{\mathfrak V}\left(\derpar{}{x^\nu}\right)=\derpar{}{x^\nu} \quad .
\eeann
Taking into account the linearity of these operators, 
these expressions allows us to compute the splitting of any vector field in $\Tan M)$.

Furthermore, for every $X\in\vf(M)$ we can obtain its horizontal lifting by doing
$({\mathfrak H}(X))(u_p)={\mathfrak h}(X(p))$, for every $u_p\in\Tan M$. 
Locally, if $\displaystyle X=X^\nu\derpar{}{x^\nu}$, we have that
$$
{\mathfrak H}(X)=
X^\nu\derpar{}{x^\nu}-\Gamma_{\nu\mu}^\rho\,X^\nu\,v^\mu\,\derpar{}{v^\rho}=
X^\nu\left(\derpar{}{x^\nu}-\Gamma_{\nu\mu}^\rho\,v^\mu\,\derpar{}{v^\rho}\right) \ .
$$

Finally, for every $X,Y\in\vf(M)$, the local expression of the covariant derivative $\nabla_XY$ is
$$
\nabla_XY=
X^\nu\derpar{Y^\rho}{x^\nu}\derpar{}{x^\rho}+\Gamma_{\nu\mu}^\rho\,X^\nu\,Y^\mu\,\derpar{}{v^\rho} \ ,
$$
and it can be expressed by means of the horizontal-vertical splitting as follows:

\begin{prop}
For every $p\in M$,
$(\nabla_XY)(p)=[\Tan_p\tau({\mathfrak h}\circ Y^C\circ X)](p)$;
where $Y^C$ denotes the {\sl complete lift} of $Y$ to $\Tan M$.
\end{prop}
\proof
Taking into account that
\beann
Y^C&=&Y^\nu\derpar{}{x^\nu}+\derpar{Y^\rho}{x^\mu}\,v^\mu\,\derpar{}{v^\rho}=
Y^\nu\left(\derpar{}{x^\nu}-\Gamma_{\nu\mu}^\rho\,v^\mu\,\derpar{}{v^\rho} \right)+
\left(\derpar{Y^\rho}{x^\mu}\,v^\mu+\Gamma_{\nu\mu}^\rho\,Y^\nu\,v^\mu\right)\derpar{}{v^\rho}
\\ &=&
{\mathfrak H}(Y^C)+{\mathfrak V}(Y^C) \quad ;
\eeann
as ${\rm h}_p^{X_p}\colon\Tan_pM\to\Tan_{X_p}(\Tan M)$
is an isomorphism,  we have that
$$
({\rm h}_p^{X_p})^{-1}([{\mathfrak V}(Y^C)]_{X_P})=
\left(\derpar{Y^\rho}{x^\mu}\,X^\mu+\Gamma_{\nu\mu}^\rho\,Y^\nu\,X^\mu\right)\derpar{}{v^\rho}\Big\vert_p.
$$
and the results follows.
\qed

\section*{Acknowledgments}

We acknowledge the financial support of the 
{\sl Ministerio de Ciencia e Innovaci\'on} (Spain), projects
MTM 2014-54855--P, MTM 2015-69124-REDT,
and of
{\sl Generalitat de Catalunya}, project 2017-SGR932.
We also acknowledge specially the help of Prof. L.A. Ibort who
introduced us in the idea of characterizing linear connections by using
the Liouville vector field.
We want to thank to all the participants in the interdisciplinary seminar on
topics of Theoretical Physics and Applied Mathematics for attending
these lectures. Their valuable suggestions and comments have allowed us
to improve these notes.

\end{document}